\documentclass[11pt]{article}
\usepackage[T1]{fontenc}
\usepackage[utf8]{inputenc}
\usepackage[english]{babel}
\usepackage{amsthm, amsmath, amsfonts, amssymb, mathtools, xcolor,
  float, hyperref, enumitem, pifont, mathrsfs, tabularx, mathabx,
  tikz, microtype, siunitx, circuitikz, microtype}
\usepackage[normalem]{ulem}
\usetikzlibrary{decorations.markings, decorations.pathmorphing,shapes,arrows,positioning}
\usepackage{comment}
\newcommand{\pp}{\mathcal{P}}
\newcommand{\g}{\mathcal{G}}

\newcommand{\zz}{\mathbb{Z}}

\renewcommand{\aa}{\mathcal{A}}

\newcommand{\cc}{\mathcal{C}}
\newcommand{\ppinv}{\overset{\leftarrow}{\pp}}
\newcommand{\pg}{\text{pre}(\g)}
\makeatletter

\definecolor{bleuclair}{RGB}{186,183,229}
\definecolor{rougeclair}{RGB}{255,230,231}
\definecolor{bleufonce}{RGB}{59,50,114}
\definecolor{rougefonce}{RGB}{127,0,4}
\newcommand\oeis[1]{\href{https://oeis.org/#1}{#1}}

\let\leq\leqslant
\let\geq\geqslant

\title{
Grand Dyck paths with air pockets
}
\author{Jean-Luc \textsc{Baril}, Sergey \textsc{Kirgizov}, R{\'e}mi \textsc{Mar{\'e}chal},\\ and Vincent \textsc{Vajnovszki}\footnote{\tt E-mails:~\{barjl,sergey.kirgizov,vvajnov\}@u-bourgogne.fr, remi.marechal01@etu.u-bourgogne.fr} \\
\rm LIB, Universit{\'e} de Bourgogne Franche-Comt{\'e}\protect\\
  B.P. 47 870, 21078 Dijon Cedex France\protect\\
}

\date{\today}

\DeclareMathAlphabet{\mymathbb}{U}{BOONDOX-ds}{m}{n}

\begin{document}

\newtheorem{defn}{Definition}
\newtheorem{prop}{Proposition}
\newtheorem{cor}{Corollary}
\newtheorem{conj}{Conjecture}
\newtheorem{rem}{Remark}
\newtheorem{lemme}{Lemma}
\newtheorem{thm}{Theorem}
\newtheorem{ex}{Example}

\maketitle

\begin{abstract}
 Grand Dyck paths with air pockets (GDAP) are a generalization of Dyck paths with air pockets by allowing them to go below the $x$-axis. We present enumerative results on GDAP (or their prefixes) subject to various restrictions such as maximal/minimal height, ordinate of the last point and particular first return decomposition. In some special cases we give bijections with other known combinatorial classes.
\end{abstract}

\section{Introduction}

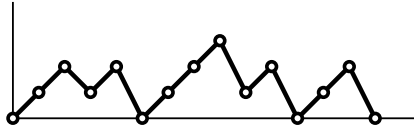
\begin{figure}[h]
\begin{center}
\scalebox{0.86}{\begin{tikzpicture}[ultra thick]
   \draw[black, thick] (0,0)--(6.2,0);\draw[black, thick] (0,0)--(0,1.8);
   \draw[black, line width=2pt](0,0)--(0.4,0.4)--(0.8,0.8)--(1.2,0.4)--(1.6,0.8)--(2,0)--(2.4,0.4)--(2.8,0.8)--(3.2,1.2)--(3.6,0.4)--(4,0.8)--(4.4,0)--(4.8,0.4)--(5.2,0.8)--(5.6,0);
  % \draw  (1.2,0.8) node {\large $C_2$};\draw  (4.8,0.8) node {\large $C_2$};

   \tikzset{every node/.style={circle, draw,fill=white,inner sep=1.5pt}}
   \node at (0,0) {};
   \node at (0.4,0.4) {};
   \node at (0.8,0.8) {};
   \node at (1.2,0.4) {};
   \node at (1.6,0.8) {};
   \node at (2,0) {};
   \node at (2.4,0.4) {};
   \node at (2.8,0.8) {};
   \node at (3.2,1.2) {};
   \node at (3.6,0.4) {};
   \node at (4,0.8) {};
   \node at (4.4,0) {};
   \node at (4.8,0.4) {};
    \node at (5.2,0.8) {};
    \node at (5.6,0) {};
\end{tikzpicture}}
\end{center}
\caption{The Dyck path with air pockets $UUDUD_2UUUD_2UD_2UUD_2$.}
\label{fig1}
\end{figure}

In a recent paper \cite{bakimava}, the authors introduce and study a new class of lattice paths, called {\it Dyck paths with air pockets} ({\it DAP} for short). 
Such a path is  a non empty lattice path in the first quadrant of $\zz^2$ starting at the origin, ending on the $x$-axis, and consisting of up-steps $U=(1,1)$ and down-steps $D_k=(1,-k)$, $k\geq 1$, where two down steps cannot be consecutive. See Figure \ref{fig1} for an example.
These paths can be viewed as ordinary Dyck paths where each maximal run of down-steps is
condensed into one large down step. As mentioned in \cite{bakimava}, they also correspond to a stack evolution with (partial) reset operations that cannot be consecutive (see for instance \cite{krin}).
In this paper, we  generalize these paths to {\it grand Dyck path with air pockets} ({\it GDAP} for short), which have the same definition as $DAP$, except that they can go below the $x$-axis, and the empty path $\varepsilon$ is considered as a GDAP.
 
 The main goal is to make enumerative studies on  GDAP (or prefix of these paths)  with various restrictions on the maximal height reached, minimal height reached, height of the last point, ...

The remaining of this paper is structured as follows.
The next section recalls some useful results from \cite{bakimava} and introduces several notations for particular subsets of GDAP. The main result of Section 3 is Theorem \ref{main_S1} which gives the  generating function counting the GDAP with respect to the length, which is the `grand' counterpart of the generating function in (\ref{rel1}) for DAP. In Section 4 we give similar results for prefixes of GDAP ending at a given ordinate,
the corresponding problem for DAP being already solved in \cite{prodinger}. In Section 5, we provide generating functions for the number of  partial GDAP that never go below the line $y=m$, with respect to the ordinate of the last point. In Section 6, we count partial GDAP lying between the lines $y=0$ and $y=t$, which correspond to partial DAP bounded by a given height $t>0$. We present a constructive bijection between these paths of length $n$ for $t=2$ and the set of compositions of $n-2$ such that no two consecutive parts have the same parity. In Section 7, we count partial GDAP lying between the lines $y=-t$ and $y=t$, and we present a constructive bijection between these paths of length $n$ for $t=1$ and the set of compositions of $n+3$ such that the first part is odd, the last part is even, and no two consecutive parts have the same parity. Finally, in Section 8  we provide enumerative results for DAP with a special first return decompostion, which proves that there are in one-to-one correspondence with Motzkin paths avoiding the patterns $UH$, $HU$ and $HH$. We leave as an open question the problem of finding a constructive bijection between these two sets.

\section{Definitions and notations}

\subsection{DAP}
 The \emph{length} of a path is the number of its steps, and for $n\geq 0$, let $\aa_n$ be the set of $n$-length DAP. By definition  $\aa_0=\aa_1=\varnothing$ and 
 we set $\aa=\bigcup_{n\geq 2}\aa_n$, see \cite{bakimava}.
A DAP is called \emph{prime} whenever it ends with $D_k$, $k\geq 2$, and
returns to the $x$-axis only once.
%For instance, the only $1$-prime Dyck path %with air pockets is $UD_1=UD$. 
The set of all prime  DAP of length $n$ is denoted $\pp_n$.
Notice that $UD$ is not prime, where for short we denote $D_1$ by $D$, so we set $\pp=\bigcup_{n\geq 3}\pp_n$.
If $\alpha=U\beta UD_k\in\pp_n$, then $2\leq k<n$  and $\beta$ is a (possibly empty) prefix of a path in $\mathcal{A}$, and we define the DAP $\alpha^\flat=\beta UD_{k-1}$, called the `lowering' of $\alpha$. For example, the path $\alpha=UUDUUD_3$ is prime, and $\alpha^\flat=UDUUD_2$. The map $\alpha\mapsto\alpha^\flat$ is clearly a bijection from $\pp_n$ to $\aa_{n-1}$ for all $n\geq 3$, and we denote by $\gamma^\sharp$ the inverse image of $\gamma\in \mathcal{A}_{n-1}$ ($\alpha^\sharp$ is a kind of `elevation' of $\alpha$).

Any DAP $\alpha\in\aa$ can be decomposed depending on its \emph{second-to-last return to the $x$-axis}: either ($i$) $\alpha=UD$, or  ($ii$) $\alpha=\beta UD$ with $\beta\in\aa$, or ($iii$) $\alpha\in \pp$, or $(iv)$  $\alpha=\beta\gamma$ where  $\beta\in\aa$ and $\gamma\in\pp$. So, if $A(x)=\sum_{n\geq 2}a_nx^n$ where $a_n$ is the cardinality of $\aa_n$, and  $P(x)=\sum_{n\geq 3}p_nx^n$ where $p_n$ is the cardinality of $\pp_n$, then we have $P(x)=xA(x)$ and the previous
decompositions imply the functional equation
$A(x)=x^2+x^2A(x)+xA(x)+xA(x)^2,$ and
\begin{equation}A(x)={\frac {1-x-{x}^{2}-\sqrt {{x}^{4}-2\,{x}^{3}-{x}^{2}-2\,x+1}}{2x
}},
\label{rel1}\end{equation}
which generates the generalized Catalan
numbers (see 
\href{https://oeis.org/A004148}{A004148} 
in~\cite{oeis}). The 
first values of $a_n$ for $2\leq n\leq 10$ 
are
$1,1,2,4,8,17,37,82,185$. In \cite{bakimava}, the authors study the enumeration of these paths according to many parameters, and they give a constructive bijection between these paths and peakless Motzkin paths (i.e. lattice paths in the first quadrant, starting at $(0,0)$, ending on the $x$-axis, made of steps $U=(1,1)$, $D=(1,-1)$ and $H=(1,0)$, and avoiding peaks of the form  $UD$).

\subsection{GDAP}
\begin{figure}[h]
\begin{center}
\scalebox{0.76}{\begin{tikzpicture}[ultra thick]
   \draw[black, thick] (0,0)--(6.2,0);\draw[black, thick] (0,-1)--(0,1.8);
   \draw[black, line width=2pt](0,0)--(0.4,0.4)--(0.8,0.8)--(1.2,0.4)--(1.6,0.8)--(2,-0.8)--(2.4,-0.4)--(2.8,-0.8)--(3.2,-0.4)--(3.6,-0)--(4,0.4)--(4.4,0.8)--(4.8,1.2)--(5.2,-0.4)--(5.6,0);
  % \draw  (1.2,0.8) node {\large $C_2$};\draw  (4.8,0.8) node {\large $C_2$};

   \tikzset{every node/.style={circle, draw, fill = white, inner sep = 1.5pt}}
   \node at (0,0) {};
   \node at (0.4,0.4) {};
   \node at (0.8,0.8) {};
   \node at (1.2,0.4) {};
   \node at (1.6,0.8) {};
   \node at (2,-0.8) {};
   \node at (2.4,-0.4) {};
   \node at (2.8,-0.8) {};
   \node at (3.2,-0.4) {};
   \node at (3.6,-0) {};
   \node at (4,0.4) {};
   \node at (4.4,0.8) {};
   \node at (4.8,1.2) {};
    \node at (5.2,-0.4) {};
    \node at (5.6,0) {};
   
\end{tikzpicture}}\qquad
\scalebox{0.76}{\begin{tikzpicture}[ultra thick]
\draw[black, thick] (0,0)--(6.2,0);\draw[black, thick] (0,-1)--(0,1.8);
   \draw[black, line width=2pt](0,0)--(0.4,0.4)--(0.8,0.8)--(1.2,0.4)--(1.6,0.8)--(2,-0.8)--(2.4,-0.4)--(2.8,-0.8)--(3.2,-0.4)--(3.6,-0)--(4,0.4)--(4.4,0.8)--(4.8,1.2)--(5.2,-0.4)--(5.6,0)--(6,0.4);
  % \draw  (1.2,0.8) node {\large $C_2$};\draw  (4.8,0.8) node {\large $C_2$};

   \tikzset{every node/.style={circle, draw, fill = white, inner sep = 1.5pt}}
   \node at (0,0) {};
   \node at (0.4,0.4) {};
   \node at (0.8,0.8) {};
   \node at (1.2,0.4) {};
   \node at (1.6,0.8) {};
   \node at (2,-0.8) {};
   \node at (2.4,-0.4) {};
   \node at (2.8,-0.8) {};
   \node at (3.2,-0.4) {};
   \node at (3.6,-0) {};
   \node at (4,0.4) {};
   \node at (4.4,0.8) {};
   \node at (4.8,1.2) {};
    \node at (5.2,-0.4) {};
    \node at (5.6,0) {};
        \node at (6,0.4) {};
  
\end{tikzpicture}}
\end{center}
\caption{A grand Dyck path with air pockets $UUDUD_4UDUUUUUD_4U$, and a partial GDAP ending at ordinate 1 with an up-step.}
\label{fig2}
\end{figure}
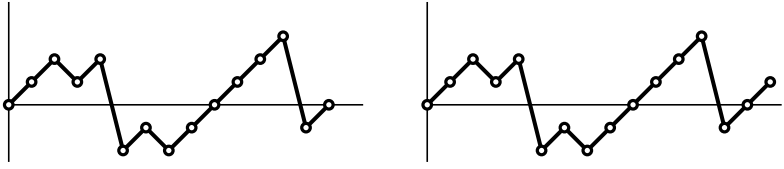

The main object of study in this paper are grand Dyck paths with air pockets (GDAP for short) which generalize DAP by allowing such paths to go below the $x$-axis; and for convenience the empty path $\varepsilon$ is a GDAP. See the first path in Figure \ref{fig2} for an example.
Let $\mathcal{G}_n$ be the set of GDAP of length $n\geq 0$, and we set $\mathcal{G}=\bigcup_{n\geq 0}\mathcal{G}_n$.

We introduce notations for several subsets of $\g$ used in this study:
\begin{itemize}
    \item $\g^+$ is the set of GDAP starting with $U$, and the empty path;
    \begin{itemize}
        \item $\g^+_1$ is the set of elements of $\g^+$ ending with $D_k, k\geq 1$;
        \item $\g^+_2$ is the set of elements of $\g^+$ ending with $U$.
    \end{itemize}
    \item $\g^-$ is the set of GDAP starting with $D_k, k\geq 1$;
    \begin{itemize}
        \item $\g^-_1$ is the set of elements of $\g^-$ ending with $D_k, k\geq 1$;
        \item $\g^-_2$ is the set of elements of $\g^-$ ending with $U$.
    \end{itemize}
\end{itemize}

Obviously, we have $\g^+=\{\varepsilon\}\cup \g^+_1\cup\g^+_2$, $\g^-= \g^-_1\cup\g^-_2$, and $\g=\g^+\cup\g^-$.
Also, we denote by $\ppinv$ the set of GDAP obtained from a prime DAP in $\pp$ by mirroring it. For instance, the mirror of $U^3D_2UD_2\in\pp$ is $D_2UD_2U^3\in\ppinv$. 

\section{Enumeration of GDAP}
\label{sec3}

In this section, we present a generating function that counts GDAP with respect to the length. 

Any element of $\g^+_1$ is of the form $\alpha\beta$, where $\alpha\in\g^+$ and $\beta\in\pp\cup\{UD\}$.
Any element of $\g^+_2$ is either of the form ($i$) $\alpha\beta$ where $\alpha\in\g^+_2$ and $\beta\in\ppinv\cup\{DU\}$, or of the form ($ii$) $\overline{\alpha\beta}$ where $\alpha\in\g^+_1$, $\beta\in\ppinv\cup\{DU\}$, and $\overline{\alpha\beta}$ is the path obtained by merging the last step of $\alpha$ together with the first step of $\beta$, i.e. if $\alpha=xD_i$ and $\beta=D_jy$, then $\overline{\alpha\beta}=xD_{i+j}y$.
Finally, any element of $\g^+$ is either the empty path $\varepsilon$, an element of $\g^+_1$, or an element of $\g^+_2$.

Thus, we deduce the following system of equations:
$$
\left\{
\begin{array}{ll}
G^+_1(x)=&G^+(x)(x^2+P(x))\\
G^+_2(x)=&G^+_2(x)(x^2+P(x))+\displaystyle\frac{1}{x}G^+_1(x)(x^2+P(x))\\
G^+(x)=&1+G^+_1(x)+G^+_2(x),
\end{array}
\right.
$$
where $P$, $G^+_1$, $G^+_2$, and $G^+$ are the generating functions with respect to the length for the cardinalities of $\pp$, $\g^+_1$, $\g^+_2$, and $\g^+$, respectively. We have $P(x)=xA(x)$, where $A$ is the generating function for the set $\mathcal{A}$ of Dyck paths with air pockets, see (\ref{rel1}). Solving the system, we get:
$$
\def\arraystretch{2.4}
\begin{array}{ll}
G^+_1(x)=&\displaystyle\frac{x^2}{\sqrt{x^4 - 2x^3 - x^2 - 2x + 1}},\\
G^+_2(x)=&\displaystyle\frac{(1-x-x^2)\sqrt{x^4 - 2x^3 - x^2 - 2x + 1} - x^4 + 2x^3 + x^2 + 2x - 1}{2x^4 - 4x^3 - 2x^2 - 4x + 2},\\
G^+(x)=&\displaystyle\frac{(1-x+x^2)\sqrt{x^4 - 2x^3 - x^2 - 2x + 1} + x^4 - 2x^3 - x^2 - 2x + 1}{2x^4 - 4x^3 - 2x^2 - 4x + 2},
\end{array}
$$
and the first terms of the respective series  expansions associated with those generating functions are:

$\bullet~ x^2 + x^3 + 2x^4 + 5x^5 + 11x^6 + 26x^7 + 63x^8 + 153x^9 + 376x^{10} + O(x^{11}),$

$\bullet~ x^3 + 2x^4 + 5x^5 + 13x^6 + 32x^7 + 80x^8 + 201x^9 + 505x^{10} + O(x^{11}),$

$\bullet~ 1 + x^2 + 2x^3 + 4x^4 + 10x^5 + 24x^6 + 58x^7 + 143x^8 + 354x^9 + 881x^{10} + O(x^{11}).$

\noindent
They correspond to the OEIS sequences \oeis{A051286}, \oeis{A110320}, and \oeis{A110236}.\\

On the other hand, any element of $\g^-$ is of the form $\alpha\beta$, where $\alpha\in\ppinv\cup\{DU\}$ and $\beta\in\g$.
Any element of $\g$ is either an element of $\g^+$ or an element of $\g^-$.

Thus, we deduce the following system of equations:
$$
\left\{
\begin{array}{rl}
G^-(x)=&(x^2+P(x))G(x)\\
G(x)=&G^+(x)+(x^2+P(x))G(x),
\end{array}
\right.
$$
where $G^-$ and $G$ are the generating functions with respect to the length for the cardinalities of  $\g^-$ and $\g$, respectively. Solving the system, we get:
$$
G^-(x)=\displaystyle\frac{(1-x+x^2 -R)(x^4  - 2x^3 - x^2 - 2x + 1 +(1-x+x^2)R)}{2(1+x-x^2 + R )R}$$ 
with $R=\sqrt{x^4 - 2x^3 - x^2 - 2x + 1}$, and we have the following result.

\begin{thm} 
\label{main_S1}
The o.g.f. that counts the set $\mathcal{G}$ with respect to the length is given by  
$$G(x)=\displaystyle\frac{x^4  - 2x^3  - x^2  - 2x + 1+(1-x+x^2)R}{(1+x-x^2 + R)R}$$
with $R$ defined as above.
\end{thm}

Notice that there is a bijection $\chi$ between the sets $\g^+_1$ and $\g^-_2$ defined as follows: for $\alpha\in \g^+_1$, $\chi(\alpha)$ is simply the mirror of $\alpha$, for instance ${\chi(UUUD_3)=D_3UUU}$.

So, we easily have
$$
\left\{
\begin{array}{rl}
G^-_2(x)=&G^+_1(x)\\
G^-_1(x)=&G^-(x)-G^-_2(x).\\
\end{array}
\right.
$$

The first terms of the series expansions  of $G^-$, $G$, $G^-_1$ and $G^-_2$ are respectively

$\bullet~x^2 + x^3 + 3x^4 + 7x^5 + 16x^6 + 39x^7 + 95x^8 + 233x^9 + 577x^{10} + O(x^{11}),$

$\bullet~1 + 2x^2 + 3x^3 + 7x^4 + 17x^5 + 40x^6 + 97x^7 + 238x^8 + 587x^9 + 1458x^{10} + O(x^{11}),$ 

$\bullet~ x^4+2x^5+5x^6+13x^7+32x^8+80x^9+201x^{10}+O(x^{11}),$ 

$\bullet~ x^2 + x^3 + 2x^4 + 5x^5 + 11x^6 + 26x^7 + 63x^8 + 153x^9 + 376x^{10} + O(x^{11}).$

\noindent
They correspond to the OEIS sequences \oeis{A203611}, \oeis{A051291}, \oeis{A110320}, and \oeis{A051286}.\\

\section{Partial GDAP ending at a given ordinate}
Let  $\pg$ be the set of partial GDAP, i.e. the set of all prefixes of elements of $\g$, see the second path in Figure \ref{fig2} for an example.
In this part, we enumerate partial GDAP ending at a given ordinate with an up-step (resp. with a down-step) with respect to the length.  Let $f_k$ (resp. $g_k$) be the generating function
for the number 
(with respect to the length)
of partial GDAP ending at ordinate $k\in\zz$ with an up-step (resp. with a down-step).  For short, we will write $f_k$ and $g_k$ instead of  $f_k(x)$ and $g_k(x)$.

According to the results in Section 2, for $k=0$ we obviously have:

$$f_0=G^+_2(x)+G^-_2(x)={\frac {1-x+{x}^{2}+\sqrt {{x}^{4}-2\,{x}^{3}-{x}^{2}-2\,x+1}}{
2\sqrt {{x}^{4}-2\,{x}^{3}-{x}^{2}-2\,x+1}}}-1
,$$
and 
$$g_0=G^+_1(x)+G^-_1(x)={\frac { \left( 1+x-{x}^{2}-\sqrt {{x}^{4}-2\,{x}^{3}-{x}^{2}-2\,x
+1}\right) x}{2\sqrt {{x}^{4}-2\,{x}^{3}-{x}^{2}-2\,x+1}}}
.$$

The first terms of the series expansions  of $f_0$ and $g_0$ are respectively

$\bullet~x^2 + 2x^3 + 4x^4 + 10x^5 + 24x^6 + 58x^7 + 143x^8 + 354x^9 + 881x^{10} + O(x^{11}),$

$\bullet~x^2 + x^3 + 3x^4 + 7x^5 + 16x^6 + 39x^7 + 95x^8 + 233x^9 + 577x^{10} + O(x^{11}).$

\noindent They correspond to the sequences
\oeis{A110236} and \oeis{A203611} in OEIS. Clearly, we have $G(x)=1+f_0+g_0$.
\medskip

For $k>0$, a partial GDAP ending at ordinate $k$  can be written $\alpha\beta$, where $\alpha$ is either empty or a GDAP ending on the $x$-axis with an up-step, and $\beta$ is a partial DAP ending at ordinate $k$. Then, we obtain
$$f_k+g_k=(1+f_0)\cdot T_k(x)$$
with $$
T_k(x)=x^ks_2^{k+1}, \mbox{ and } s_2=\frac{1+x-x^2-\sqrt{-x^2-2x^3-2x+x^4+1}}{2x},
$$ where $T_k$ is the o.g.f. that counts DAP ending at ordinate $k$ with respect to length, which is already obtained in \cite{prodinger}.

As a consequence, we deduce the following result.

\begin{thm} The o.g.f. that counts the partial GDAP ending at a positive ordinate (with respect to the length) is given by 
$${\frac { \left( {x}^{2}-x-1+\sqrt {{x}^{4}-2\,{x}^{3}-{x}^{2}-2\,x+
1} \right) ^{2}}{4x\sqrt {{x}^{4}-2\,{x}^{3}-{x}^{2}-2\,x+1}}}
.$$
\end{thm}

\noindent
The first terms of the series expansion are:
$x+x^2+4x^3+9x^4+22x^5$+$55x^6+136x^7+339x^8+849x^9+2132x^{10}+ O(x^{11})$.

\medskip

For $k<0$, a partial GDAP ending at ordinate $k$ can be written $\beta$ or $\overline{\alpha\beta}$,
where $\alpha$ is a GDAP ending with a down-step, and $\beta$ is the symmetric about the $x$-axis of a partial DAP ending at ordinate $-k>0$ in the right-to-left  model  studied in \cite{prodinger}. Then, we obtain $$f_k+g_k=R_k(x)\left(1+\frac{g_0(x)}{x}\right),$$ with 
$$
R_k(x)=(s_2-1)\cdot\frac{s_2^{-k-1}}{x}.
$$
For instance, the first terms of the series expansion for $k=-1,-2$ are

 $\bullet~x+2x^2+4x^3+10x^4+24x^5+58x^6+143x^7+354x^8+881x^9+2204x^{10}+O(x^{11}),$ 

$\bullet~x+2x^2+5x^3+13x^4+32x^5+80x^6+201x^7+505x^8+1273x^9+3217x^{10}+O(x^{11})$  
 which correspond to the sequences \oeis{A110236} and \oeis{A110320} in OEIS.

Notice  that partial GDAP of length $n-1$ and ending at ordinate $k=-1$ are in one-to-one correspondence with non-empty GDAP of length $n$ and starting
with an up-step (see the set $\mathcal{G}^+$ from
Section~\ref{sec3}).  To perceive it, one can add an up-step at the end of a partial GDAP ending at $k=-1$, apply a symmetry about the $x$-axis, and consider the mirror of this path.

  Moreover, $x^2 \cdot G^+_2(x)$ (see Section~\ref{sec3}) is the
   generating function for partial GDAP ending at $k=-2$.  Such
   partial GDAP of length $n$ are in one-to-one correspondence  with paths of length $n+2$  in $\mathcal{G}^+_2$. To see it,  from a partial GDAP ending at $k=-2$, one can add an up-step at the beginning of the path and another one at the end.

Since there is an infinite number
of partial paths of length $n$
ending at negative height, 
we cannot provide an ordinary generating function (with respect to the length)
for these paths. So, we get around this in the next section by counting partial GDAP
lying above the line $y=m$ for a given $m\leq0$.

%%%%%%%%%%
\begin{comment}
Now, we introduce the bivariate generating functions
$$
F(u,z)=F(u)=\sum_{k=-\infty}^\infty u^kf_k \quad\text{and}\quad G(u,z)=G(u)=\sum_{k=-\infty}^\infty u^kg_k.
$$
Making use of the recursions above, we get:
$$
\begin{array}{rcl}
F(u)&=&1+zu\left(F(u)+G(u)\right),\\
G(u)&=&\displaystyle\frac{z}{u-1}F(u).
\end{array}
$$
Plugging the second equation into the first one, we get:
$$
F(u)=\frac{1-u}{zu^2+(z^2-z-1)u+1}.
$$

{\color{red}{\bf ça me donne des séries avec des signes $+$ et des signes $-$ pour $F$, $G$, et $F+G$}}
\end{comment}
%%%%%%%%%%

%tout ça ne marche pas car calcul illicite (module d'une série géométrique plus grand que 1)

\section{Minorized partial GDAP}

Let us denote by $\pg_m$ the set of partial GDAP which never go below the line  $y=m$, $m\leq 0$, and let us reuse the same notations as in the previous section for the generating functions $f_k$ and $g_k$ in this subset of $\pg$. Obviously, we have $f_k=0$ for all $k\leq m$ and $g_k=0$ for all $k<m$. By convenience, we count the empty path in $f_0$. Then, the o.g.f.'s satisfy the following equations:
$$
\left\{
\begin{array}{rcl}
f_0&=&1+xf_{-1}+xg_{-1},\\
\forall k\geq m+1, k\neq 0,\quad f_k&=&xf_{k-1}+xg_{k-1},\\
\forall k\geq m,\quad g_k&=&\sum_{i=1}^\infty xf_{k+i}.
\end{array}
\right.
$$
As a consequence, we have $f_{m+1}=xg_m$. Now, we introduce the bivariate generating functions
$$
f(u,x)=f(u)=\sum_{k=m+1}^\infty u^kf_k \quad\text{and}\quad g(u,x)=g(u)=\sum_{k=m}^\infty u^kg_k.
$$
Making use of the recursions above, we get:
$$
\begin{array}{rcl}
f(u)&=&1+xu\left(f(u)+g(u)\right),\\
g(u)&=&\displaystyle\frac{x}{1-u}\left(u^mf(1)-f(u)\right).
\end{array}
$$
Plugging the second equation into the first one, we get:
$$
f(u)=1+xuf(u)+\displaystyle\frac{x^2u}{1-u}\left(u^mf(1)-f(u)\right).
$$
Solving for $f(u)$, we finally get:
\begin{equation}
f(u)=\displaystyle\frac{1-u+x^2u^{m+1}f(1)}{1-u-xu+xu^2+x^2u}\label{fu}.
\end{equation}

In order to compute $f(1)$, we use the kernel method (see~\cite{ban, pro}) on $f(u)$. We can rewrite the denominator---which is a polynomial in $u$, of degree $2$---as $x(u-r_1)(u-r_2)$, where:
$$
r_1=\frac{1+x-x^2+\sqrt{x^4-2x^3-x^2-2x+1}}{2x},
$$

$$r_2=\frac{1+x-x^2-\sqrt{x^4-2x^3-x^2-2x+1}}{2x},$$

%%%%%%%%%%
\begin{comment}
Using one of Vieta's formulas, we have:
$$
r+\bar{r}=1-x+\frac{1}{x}.
$$

Hence:
$$
x(u-r)=x\bar{r}-x+x^2-1+xu, %ne servira que lorsqu'on aura une expression du numérateur divisé par u-\bar{r}
$$
\end{comment}
%%%%%%%%%%

\noindent
and then, relation~\eqref{fu} implies 
$$
f(u)\cdot(x(u-r_1)(u-r_2))=1-u+x^2u^{m+1}f(1).
$$
Plugging $u=r_2$ (which has a Taylor expansion at $x=0$), we obtain:
$$
1-r_2+x^2r_2^{m+1}f(1)=0,
$$
which gives an expression for $f(1)$:
$$
f(1)=\frac{r_2-1}{x^2r_2^{m+1}},
$$
and then:
$$
\def\arraystretch{2.4}
\begin{array}{c}
f(u)=\displaystyle\frac{1-u+(r_2-1)\left(\frac{u}{r_2}\right)^{m+1}}{1-u-xu+xu^2+x^2u},\\
g(u)=\displaystyle\frac{x}{1-u}\left(u^mf(1)-f(u)\right).
\end{array}$$

%This means we can simplify both the numerator and the denominator in $F(u)$ by the factor $u-\bar{r}$.
Finally, we have:

\begin{thm} The o.g.f. that counts the partial GDAP above the line $y=m$ (with respect to the length) is given by 
$$f(1)+g(1)=\frac{r_2^{-m}-r_2^{-1-m}-x^2}{x^3}.$$
\end{thm}

\noindent
For instance, if $m=-1,-2$ the first terms of the series expansions are

$\bullet$ $1+2x+4x^2+8x^3+17x^4+37x^5+82x^6+185x^7+423x^8+978x^9+2283x^{10}+O(x^{11})$,

$\bullet$ $1+3x+6x^2+13x^3+29x^4+65x^5+148x^6+341x^7+793x^8+1860x^9+4395x^{10}+O(x^{11})$.

\noindent
They correspond to the sequences \oeis{A004148} and \oeis{A093128} in OEIS.

\section{Partial (G)DAP bounded by $y=0$ and $y=t$}

\subsection{Enumerative results}
\label{enumeration_0_t}

In this section, we count partial GDAP lying between the lines $y=0$ and  $y=t$, which correspond to partial DAP bounded by a given height $t>0$. We introduce the notation $f_k^t$, $g_k^t$ for $0\leq k \leq t$, $f^t(u)$, and $g^t(u)$, which are the counterparts of $f_k$, $g_k$, $f(u)$, and $g(u)$ defined in the previous section. So, we deduce the following system of equations:
$$
\left[\begin{array}{cccc|cccc}
-1&&&&0&&&\\
x&-1&&&x&0&&\\
&\ddots&\ddots&&&\ddots&\ddots&\\
&&x&-1&&&x&0\\
\hline
0&x&\hdots&x&-1&&&\\
&0&\ddots&\vdots&&\ddots&&\\
&&\ddots&x&&&\ddots&\\
&&&0&&&&-1
\end{array}\right]
\cdot
\begin{bmatrix}
f_0^t\\
\vdots\\
\vdots\\
f_t^t\\
\hline
g_0^t\\
\vdots\\
\vdots\\
g_t^t
\end{bmatrix}
=
\begin{bmatrix}
-1\\
0\\
\vdots\\
\vdots\\
\vdots\\
\vdots\\
\vdots\\
0
\end{bmatrix}.
$$
For a given height $t\geq 0$, the previous matrix (denoted $A_t$) is square with $2(t+1)$ rows. Using classical properties of the determinant (in particular $\mbox{det}\left(\begin{array}{cc}A & B\\C & D\\\end{array}\right)=\mbox{det}(AD-BC)$ whenever $D$ is invertible, and $C$ and $D$ commute \cite{Syl}), we can easily prove that $D_t=\det(A_t)$ satisfies
$$
D_{t+2}+(x^2-x-1)D_{t+1}+xD_t=0,
$$
anchored with $D_0=1$, and $D_1=1-x^2$. Then we deduce
$$
D_t=\frac{2^tx^{t+1}}{W}\left(\frac{ W-x^{2}+x -1}{\left(W-x^{2}+x +1\right)^{t+1}}
+(-1)^{t+1}\frac{W+x^{2}-x +1}{\left(W+x^{2}-x -1\right)^{t+1}}\right)
,
$$
where $$W=\sqrt{x^{4}-2 x^{3}-x^{2}-2 x +1}.$$ For instance, we have $D_2=x^{4}-x^{3}-2 x^{2}+1$, and $D_3=-x^{6}+2 x^{5}+2 x^{4}-2 x^{3}-3 x^{2}+1$.\par 
Using Cramer's rule to solve the system, for $0\leq k \leq t$, we have
\begin{equation}
f_k^t=\frac{N_k^t}{D_t},\quad g_k^t=\frac{N_{t+1+k}^t}{D_t}\label{eqnkt},
\end{equation}
where $N_k^t$ is the determinant of the matrix $A_t(k)$ obtained from $A_t$ by replacing the $(k+1)$-th column with the vector $(-1,0,\hdots,0)^T$.\par

As we have done for $D_t$, it is easy to prove that $N_k^t$ satisfies the following recurrence relations, for $0\leq k \leq t$:
$$
\left\{\begin{array}{rcll}
N_0^t&=&D_t&\\
N_{2t+1}^t&=&0&\\
N_k^t&=&xN_{k-1}^{t-1}&1\leq k\leq t\\
N_{t+k}^t&=&xN_{t+k-2}^{t-1}&2\leq k\leq t\\
N_{t+1}^t&=&x^2N_0^{t-1}+xN_t^{t-1}.&
\end{array}\right.
$$
See Table~\ref{nkt} for exact values of $N_k^t$ when $0\leq t\leq3$ and $0\leq k\leq7$.

\begin{table}[H]
    \centering
   \scalebox{1}{ \begin{tabular}{ccccc}
    $k\backslash t$&$0$&$1$&$2$&$3$\\
    \hline
$0$&
$1$&
$-x^{2}+1$&
$x^{4}-x^{3}-2 x^{2}+1$&
$-x^{6}+2 x^{5}+2 x^{4}-2 x^{3}-3 x^{2}+1$\\
%&$x^{8}-3 x^{7}-x^{6}+5 x^{5}+4 x^{4}-3 x^{3}-4 x^{2}+1$\\
$1$&
$0$&
$x$&
$-x^{3}+x$&
$x^{5}-x^{4}-2 x^{3}+x$\\
%&$-x^{7}+2 x^{6}+2 x^{5}-2 x^{4}-3 x^{3}+x$\\
$2$&&
$x^{2}$&
$x^{2}$&
$-x^{4}+x^{2}$\\
%&$x^{6}-x^{5}-2 x^{4}+x^{2}$\\
$3$&&
$0$&
$-x^{4}+x^{3}+x^{2}$&
$x^{3}$\\
%&$-x^{5}+x^{3}$\\
$4$&&&
$x^{3}$&
$x^{6}-2 x^{5}-x^{4}+x^{3}+x^{2}$\\
%&$x^{4}$\\
$5$&&&
$0$&
$-x^{5}+x^{4}+x^{3}$\\
%\\&$-x^{8}+3 x^{7}-3 x^{5}-2 x^{4}+x^{3}+x^{2}$\\
$6$&&&&
$x^{4}$\\
%&$x^{7}-2 x^{6}-x^{5}+x^{4}+x^{3}$\\
$7$&&&&
$0$\\
%&$-x^{6}+x^{5}+x^{4}$\\
    \end{tabular}}
    \caption{The first values of $N_k^t$ for $0\leq t\leq3$ and $0\leq k\leq 7$.}
    \label{nkt}
\end{table}

Using \eqref{eqnkt} and the above recurrence relations for $N_k^t$, we can deduce closed forms for $f_k^t$, $g_k^t$, $0\leq k\leq t$. 

So, we can state the following result.

\begin{thm} The o.g.f. that counts the  nonempty GDAP bounded by the lines $y=0$ and $y=t$ (with respect to the length) is
 $$g_0^t=\frac{N_{t+1}^t}{D_t}$$
with  
  $$
 N_{t+1}^t =\frac{2^{t+2}x^{t+3}(-1)^t}{W(x^2-x-1)^2-W^3}\left( \frac{1}{(x^{2}-x -1+W)^{t}}
  - \frac{1}{(x^{2}-x -1-W)^t}
\right).$$
\end{thm}

For instance, if $t=1,2,3,4$, then we have  

$$
g_0^1 = \frac{x^{2}}{1-x^{2}}, 
g_0^2 = \frac{x^{2} \left(1+x-x^{2}\right)}{x^{4}-x^{3}-2 x^{2}+1},
g_0^3 = \frac{x^{2} \left(x^{4}-2 x^{3}-x^{2}+x +1\right)}{\left(x^{3}-2 x^{2}-x +1\right) \left(1+x-x^{3}\right)},\text{ and}
$$

$$
g_0^4 = \frac{-x^{8}+3 x^{7}-3 x^{5}-2 x^{4}+x^{3}+x^{2}}{x^{8}-3 x^{7}-x^{6}+5 x^{5}+4 x^{4}-3 x^{3}-4 x^{2}+1},
$$
and the first terms of the series expansion of these generating functions are respectively

    $\bullet$~$x^{2}+x^{4}+x^{6}+x^{8}+x^{10}+O(x^{11})$,

 $\bullet$~$x^{2}+x^{3}+x^{4}+3 x^{5}+2 x^{6}+6 x^{7}+6 x^{8}+11 x^{9}+16 x^{10}+O(x^{11})$,

 $\bullet$~$x^{2}+x^{3}+2 x^{4}+3 x^{5}+7 x^{6}+9 x^{7}+22 x^{8}+32 x^{9}+66 x^{10}+O(x^{11})$,
 
 $\bullet$~$x^{2}+x^{3}+2 x^{4}+4 x^{5}+7 x^{6}+16 x^{7}+27 x^{8}+63 x^{9}+112 x^{10}+O(x^{11})$.

\noindent
The first two correspond to shifts of \oeis{A000035} and \oeis{A062200}. The last two sequences do not appear in \cite{oeis}.

\subsection{Bijection with a set of compositions}

As stated above, the enumeration of the set $\g_n^{[0,2]}$  of GDAP bounded by $y=0$ and $y=2$ is given by $g_0^2=\frac{x^{2} \left(1+x-x^{2}\right)}{x^{4}-x^{3}-2 x^{2}+1}$ which have  a series expansion where the coefficients coincide (up to a shift) with the sequence \oeis{A062200} in \cite{oeis}. In this part, for any $n\geq 0$, we exhibit a constructive bijection $\psi$ between $\g_n^{[0,2]}$ and the set  $\cc(n-2)$  of compositions of  $n-2$ such that no two consecutive parts have the same parity (see \cite{mans} for the enumeration of these objects, and \cite{mansour} for more results about the enumeration of compositions with regard to several statistics on parts).

Let us define the map $\psi$. Assuming $n\geq 2$, let  $\alpha=\alpha_1\hdots \alpha_n\in\g_n^{[0,2]}$ and $\alpha'=\alpha_2\alpha_3\hdots \alpha_{n-1}$. We write $\alpha'=B_1B_2\hdots B_r$ where each $B_i$ is a subpath of $\alpha'$ satisfying the following rules:
\begin{itemize}
     \item if $\alpha'$ does not contain  $U^2$ and $UD_2$, then $r=1$ and $B_1=\alpha'$;
      \item otherwise,  we split $\alpha'$ into 
      subpaths $B_i$, $1\leq i\leq r$, by cutting it 
      after all up-steps that are followed by another up-step or a $D_2$-step.
   \end{itemize}

For instance, if $\alpha=UUD_2UUDUD_2UDUDUUD_2=U\alpha' D_2$, then $\alpha'=B_1B_2B_3B_4B_5$ where $B_1=U$, $B_2=D_2U$, $B_3=UDU$, $B_4=D_2UDUDU$, $B_5=U$. We refer to Figure \ref{exemple1_psi} for an illustration of this decompostion.

%Note that the number $r$ of subpaths in the decomposition of $\alpha'$ is at least one %(e.g. $\alpha=(UD)^{n/2}$), and at most $2n/3$ (e.g. $\alpha=(UD_2U)^{n/3}$). 

\begin{figure}[h]
\begin{center}
\scalebox{0.8}{\begin{tikzpicture}[ultra thick]
    \fill[bleuclair] (0,0) rectangle (0.4,0.8);
    \fill[rougeclair] (0.4,0.8) rectangle (1.2,-0);
    \fill[bleuclair] (1.2,0) rectangle (2.4,0.8);
    \fill[rougeclair] (2.4,0.8) rectangle (4.8,0);
    \fill[bleuclair] (4.8,0) rectangle (5.2,0.8);
    \draw[black, thick] (-0.4,0)--(6.2,0);\draw[black, thick] (-0.4,-0.4)--(-0.4,0.8);
   \draw[black, line width=2pt](-0.4,0) --(0,0.4)--(0.4,0.8)--(0.8,0)--(1.2,0.4)--(1.6,0.8)--(2,0.4)--(2.4,0.8)--(2.8,0)--(3.2,0.4)--(3.6,0)--(4,0.4)--(4.4,0)--(4.8,0.4)--(5.2,0.8)--(5.6,0);
  % \draw  (1.2,0.8) node {\large $C_2$};\draw  (4.8,0.8) node {\large $C_2$};

\node at (0.2,1.2) {{\color{bleufonce}{$1$}}};
\node at (0.8,1.2) {{\color{rougefonce}{$2$}}};
\node at (1.8,1.2) {{\color{bleufonce}{$3$}}};
\node at (3.6,1.2) {{\color{rougefonce}{$6$}}};
\node at (5.2,1.2) {{\color{bleufonce}{$1$}}};

\draw[-to] (6.7,0)--(7.7,0);

\node at (10,0) {$1+2+3+6+1=13$};

   \tikzset{every node/.style={circle, draw,fill=white,inner sep=1.5pt}}
    \node at (-0.4,0) {};
   \node at (0,0.4) {};
   \node at (0.4,0.8) {};
   \node at (0.8,0) {};
   \node at (1.2,0.4) {};
   \node at (1.6,0.8) {};
   \node at (2,0.4) {};
   \node at (2.4,0.8) {};
   \node at (2.8,0) {};
   \node at (3.2,0.4) {};
   \node at (3.6,0) {};
   \node at (4,0.4) {};
   \node at (4.4,0) {};
   \node at (4.8,0.4) {};
    \node at (5.2,0.8) {};
    \node at (5.6,0) {};

\end{tikzpicture}}
\end{center}
\caption{The image by $\psi$ of $\alpha=UUD_2UUDUD_2UDUDUUD_2$ is $\psi(\alpha)=1,2,3,6,1$.}
\label{exemple1_psi}
\end{figure}

Let $b_1,\hdots,b_r$ be the lengths of the subpaths $B_1,\hdots, B_r$, respectively. It is clear that $b_1+b_2+\hdots+b_r=n-2$. Moreover,  if the subpath $B_i$ starts with $U$ and ends with $U$, then $B_i$ is of he form $U(DU)^k$ for some $k\geq 0$, and $B_{i+1}$ is necessarily of the form $D_2(UD)^\ell U$ for some $\ell$, which implies that $b_i$ is odd and $b_{i+1}$ is even; if the subpath $B_i$ starts with $D_2$ and ends with $U$, then $B_i$ is of the form $D_2(UD)^k U$ for some $k\geq 0$, and $B_{i+1}$ is necessarily of the form $U(DU)^\ell$ for some $\ell$, which implies that $b_i$ is even and $b_{i+1}$ is odd. Thus, two consecutive  $b_i$ and $b_{i+1}$ always have different parities. Then, the above procedure defines a map $\psi$ from $\g_n^{[0,2]}$ to the set $\mathcal{C}(n-2)$, and 
for $\alpha \in \g_n^{[0,2]}$, we set
$\psi(\alpha) = b_1, b_2, \ldots, b_r$.

\begin{thm} The map $\psi$ from $\g_n^{[0,2]}$  to  $\cc(n-2)$ is a bijection.
\end{thm}

\begin{proof}
  Since $\g_n^{[0,2]}$ and $\cc(n-2)$ have the same cardinality (see the o.g.f. $g_0^2$ and the sequence \oeis{A062200} at the end of subsection~\ref{enumeration_0_t},
  it suffices to prove that $\psi$ is surjective.

 Let $c=c_1,\hdots,c_r$ be a composition in $\cc(n-2)$, with $r\geq 2$ (the case $r=1$ being  trivial since if $c_1$ is even, then we have $c=\psi(U(DU)^{c_1/2}D)$, and if $c_1$ is odd we have $c=\psi(U(UD)^{(c_1-1)/2}UD_2)$). 
 
 For $r\geq 2$, we distinguish four cases: ($i$) $c_1$ and $c_r$ are even,  ($ii$) $c_1$ is even and $c_r$ is odd, ($iii$)  $c_1$ and $c_r$ are odd,  ($iv$) $c_1$ is odd and $c_r$ is even. According to each case, we define $\alpha\in\g_n^{[0,2]}$ such that $\psi(\alpha)=c$ as follows:

\noindent Case ($i$):
$$\alpha=
U\colorbox{rougeclair}{$(DU)^{c_1/2}$}
\colorbox{bleuclair}{$(UD)^{(c_{2}-1)/2}U$}\colorbox{rougeclair}{$D_2U(DU)^{(c_3-2)/2}$}\hdots\colorbox{rougeclair}{$D_2U(DU)^{(c_r-2)/2}$}D;
$$
\noindent Case ($ii$): 
$$\alpha=
U\colorbox{rougeclair}{$(DU)^{c_1/2}$}
\colorbox{bleuclair}{$(UD)^{(c_{2}-1)/2}U$}\colorbox{rougeclair}{$D_2U(DU)^{(c_3-2)/2}$}\hdots\colorbox{bleuclair}{$(UD)^{(c_{r}-1)/2}U$}D_2;
$$
\noindent Case ($iii$):
$$\alpha=U
\colorbox{bleuclair}{$(UD)^{(c_{1}-1)/2}U$}\colorbox{rougeclair}{$D_2U(DU)^{(c_2-2)/2}$}
\colorbox{bleuclair}{$(UD)^{(c_{3}-1)/2}U$}\hdots\colorbox{bleuclair}{$(UD)^{(c_{r}-1)/2}U$}D_2;
$$
\noindent Case ($iv$): 
$$\alpha=U
\colorbox{bleuclair}{$(UD)^{(c_{1}-1)/2}U$}\colorbox{rougeclair}{$D_2U(DU)^{(c_2-2)/2}$}
\colorbox{bleuclair}{$(UD)^{(c_{3}-1)/2}U$}\hdots\colorbox{rougeclair}{$D_2U(DU)^{(c_r-2)/2}$}D.
$$
For each case, it is clear that $\alpha$ belongs to $\g_n^{[0,2]}$, which implies that $\psi$ is surjective, and then bijective.
\end{proof}

\section{GDAP bounded by $y=-t$ and $y=t$}

\subsection{Enumerative results}\label{enumeration_-t_t}
In this section, we count GDAP lying between the lines $y=-t$ and  $y=t$. We introduce the notation $f_k^t$, $g_k^t$ for $-t\leq k \leq t$, $f^t(u)$, and $g^t(u)$, which are the counterparts of $f_k$, $g_k$, $f(u)$, and $g(u)$. So, we deduce the following system of equations:
$$
\left[\begin{array}{cccc|cccc}
-1&&&&0&&&\\
x&-1&&&x&0&&\\
&\ddots&\ddots&&&\ddots&\ddots&\\
&&x&-1&&&x&0\\
\hline
0&x&\hdots&x&-1&&&\\
&0&\ddots&\vdots&&\ddots&&\\
&&\ddots&x&&&\ddots&\\
&&&0&&&&-1
\end{array}\right]
\cdot
\begin{bmatrix}
f_{-t}^t\\
\vdots\\
f_{-1}^t\\
f_0^t\\
f_1^t\\
\vdots\\

f_t^t\\
\hline
g_{-t}^t\\
\vdots\\
\vdots\\
g_t^t
\end{bmatrix}
=
\begin{bmatrix}
0\\
\vdots\\
0\\
-1\\
0\\
\vdots\\
0\\
\hline
0\\
\vdots\\
\vdots\\
0
\end{bmatrix}.
$$
For a given height $t\geq 0$, the previous matrix (denoted $A_t^{'}$) is square with $2(2t+1)$ rows. We notice that for all $t\geq 0$, the matrix $A_t^{'}$ is identical to the matrix $A_{2t}$ defined in the previous section. Hence, we have $D_t^{'}:=\det(A_t^{'})=\det(A_{2t})=D_{2t}$, i.e.
$$
D_t^{'}=\frac{4^tx^{2t+1}}{W}\left(\frac{ W-x^{2}+x -1}{\left(W-x^{2}+x +1\right)^{2t+1}}
-\frac{W+x^{2}-x +1}{\left(W+x^{2}-x -1\right)^{2t+1}}\right)
,
$$
where $$\textcolor{red}{W=}\sqrt{x^{4}-2 x^{3}-x^{2}-2 x +1}.$$\par 
Using Cramer's rule to solve the system, for $-t\leq k \leq t$, we have
\begin{equation}
f_k^t=\frac{\widetilde{N}_k^t}{D^{'}_t},\quad g_k^t=\frac{\widetilde{N}_{2t+1+k}^t}{D^{'}_t}\label{eqnkt_bis},
\end{equation}
where $\widetilde{N}_k^t$ is the determinant of the matrix $A^{'}_t(k)$ obtained from $A^{'}_t$ by replacing the $(k+t+1)$-th column with the vector $(0,\hdots,0,-1,0,\hdots,0)^T$, where the $-1$ is in the $(t+1)$-th position.\par

 Now, we focus on the calculation of $f_k^t$ and $g_k^t$ for $k=0$. The other cases can be obtained similarly, but they are much more technical and less interesting to present them here.
With the same arguments as in the previous section (in particular, using the mentioned property of the determinant on  blocks), it is easy  to prove that $\widetilde{N}_0^t$ satisfies: 
$$\widetilde{N}_0^t=D_{t-1}\cdot D_t.$$
Moreover, we have $$\widetilde{N}_{2t+1}^t=D_{t-1}\cdot N_{t+1}^t.$$

%See Table~\ref{nkt_bis} for exact values of $\widetilde{N}_k^t$ when $0\leq t\leq3$ and $-3\leq k\leq10$.

\begin{comment}
\begin{table}[H]
    \centering
   \scalebox{0.65}{ \begin{tabular}{ccccc}
    $k\backslash t$&$0$&$1$&$2$&$3$\\
    \hline
$-3$&
&
&
&
$0$\\
$-2$&
&
&
$0$&
$-x^{8}+3 x^{7}-3 x^{5}-2 x^{4}+x^{3}+x^{2}$\\
$-1$&
&
$0$&
$x^{6}-2 x^{5}-x^{4}+x^{3}+x^{2}$&
$x^{10}-4 x^{9}+2 x^{8}+6 x^{7}-x^{6}-6 x^{5}-2 x^{4}+2 x^{3}+x^{2}$\\
$0$&
$1$&
$-x^{2}+1$&
$-x^{6}+x^{5}+3 x^{4}-x^{3}-3 x^{2}+1$&
$-x^{10}+3 x^{9}+2 x^{8}-8 x^{7}-6 x^{6}+9 x^{5}+9 x^{4}-3 x^{3}-5 x^{2}+1$\\
$1$&
$0$&
$x$&
$x^{5}-2 x^{3}+x$&
$x^{9}-2 x^{8}-3 x^{7}+4 x^{6}+6 x^{5}-2 x^{4}-4 x^{3}+x$\\
$2$&
&
$-x^{3}+x^{2}+x$&
$-x^{4}+x^{2}$&
$-x^{8}+x^{7}+3 x^{6}-x^{5}-3 x^{4}+x^{2}$\\
$3$&
&
$x^{2}$&
$x^{5}-2 x^{4}-x^{3}+x^{2}+x$&
$x^{7}-x^{6}-2 x^{5}+x^{3}$\\
$4$&
&
$0$&
$-x^{7}+2 x^{6}+2 x^{5}-3 x^{4}-2 x^{3}+x^{2}+x$&
$-x^{7}+3 x^{6}-3 x^{4}-2 x^{3}+x^{2}+x$\\
$5$&
&
&
$x^{6}-x^{5}-2 x^{4}+x^{3}+x^{2}$&
$x^{9}-3 x^{8}-x^{7}+6 x^{6}+2 x^{5}-4 x^{4}-3 x^{3}+x^{2}+x
$\\
$6$&
&
&
$-x^{5}+x^{3}$&
$-x^{11}+4 x^{10}-x^{9}-9 x^{8}+12 x^{6}+4 x^{5}-6 x^{4}-4 x^{3}+x^{2}+x$\\
$7$&
&
&
$0$&
$x^{10}-3 x^{9}-x^{8}+6 x^{7}+3 x^{6}-5 x^{5}-3 x^{4}+x^{3}+x^{2}$\\
$8$&
&
&
&
$-x^{9}+2 x^{8}+2 x^{7}-3 x^{6}-3 x^{5}+x^{4}+x^{3}$\\
$9$&
&
&
&
$x^{8}-x^{7}-2 x^{6}+x^{4}$\\
$10$&
&
&
&
$0$\\

    \end{tabular}}
    \caption{The first values of $\widetilde{N}_k^t$ for $0\leq t\leq3$ and $-3\leq k\leq10$.}
    \label{nkt_bis}
\end{table}
\end{comment}

Using the results obtained in the previous section, these two relations allow to obtain a close form for $\widetilde{N}_0^t$ and $\widetilde{N}_{2t+1}^t$. Using \eqref{eqnkt_bis}, we deduce close forms for $f_0^t$ and $g_0^t$ and we can state the following result.
\begin{thm} The o.g.f. that counts the nonempty GDAP bounded by the lines $y=-t$ and $y=t$ (with respect to the length) is
 $$f_0^t+g_0^t=\frac{D_{t-1}}{D_{2t}}\cdot\left(D_t+N_{t+1}^t\right)$$ where $D_t$ and $N_{t+1}^t$ are defined in the previous section.
\end{thm}

For instance, if $t=1,2,3$, then we have

$$
f_0^1+g_0^1 = \frac{\scriptstyle 1}{\scriptstyle x^{4}-x^{3}-2 x^{2}+1}, 
f_0^2+g_0^2 = \frac{\scriptstyle\left(x -1\right) \left(x +1\right)^{2}}{\scriptstyle x^{7}-2 x^{6}-3 x^{5}+2 x^{4}+6 x^{3}+3 x^{2}-x -1},
$$

$$
f_0^3+g_0^3 = \frac{\scriptstyle{\left(x^{4}-x^{3}-2 x^{2}+1\right)^{2}}}{\scriptstyle{x^{12}-5 x^{11}+4 x^{10}+10 x^{9}-4 x^{8}-19 x^{7}-4 x^{6}+17 x^{5}+11 x^{4}-5 x^{3}-6 x^{2}+1}},
$$

%$$
%f_0^4+g_0^4 = \frac{\left(x^{3}-2 x^{2}-x +1\right)^{2} \left(x^{3}-x -1\right)^{2}}{x^{16}-7 %x^{15}+13 x^{14}+7 x^{13}-27 x^{12}-26 x^{11}+41 x^{10}+57 x^{9}-22 x^{8}-68 x^{7}-15 x^{6}+37 %x^{5}+22 x^{4}-7 x^{3}-8 x^{2}+1},
%$$

\noindent
and the first terms of the series expansion of these generating functions are respectively

    $\bullet$~$1+2 x^{2}+x^{3}+3 x^{4}+4 x^{5}+5 x^{6}+10 x^{7}+11 x^{8}+21 x^{9}+27 x^{10}+O(x^{11})$,

 $\bullet$~$1+2 x^{2}+3 x^{3}+5 x^{4}+13 x^{5}+22 x^{6}+48 x^{7}+93 x^{8}+190 x^{9}+375 x^{10}+O(x^{11})$,

 $\bullet$~$1+2 x^{2}+3 x^{3}+7 x^{4}+15 x^{5}+36 x^{6}+75 x^{7}+176 x^{8}+386 x^{9}+869 x^{10}+O(x^{11})$.
 
% $\bullet$~$1+2 x^{2}+3 x^{3}+7 x^{4}+17 x^{5}+38 x^{6}+93 x^{7}+212 x^{8}+513 x^{9}+1194 x^{10}+O(x^{11})$,

\noindent
The first one corresponds to a shift of \oeis{A122514}. The last two sequences do not appear in \cite{oeis}.

\subsection{Bijection with a set of  compositions}

As stated above, the enumeration of the set $\g_n^{[-1,1]}$ of  GDAP bounded by $y=-1$ and $y=1$ is given by  $f_0^1+g_0^1=\frac{1}{x^{4}-x^{3}-2 x^{2}+1}$ which have a series expansion where the coefficients coincide (up to a shift) with the sequence \oeis{A122514}. In this part, for any $n\geq 0$, we exhibit a constructive bijection $\phi$ between $\g_n^{[-1,1]}$ and the set  $\cc'(n+3)$  of compositions of  $n+3$ such that the first part is odd, the last part is even, and no two consecutive parts have the same parity.

Now, let us define the map $\phi$. Assuming $n\geq 2$, let $\alpha=\alpha_1\hdots \alpha_n\in\g_n^{[-1,1]}$. We write $\alpha=B_1B_2\hdots B_r$ where each $B_i$ is a subpath of $\alpha$ obtained by applying the same decomposition made on $\alpha'$ in subsection~6.2. Let $b_1,b_2,\hdots,b_r$ be the lengths of subpaths $B_1,\hdots, B_r$ respectively. 
In the case $r\geq 2$, let $L$ be the reversed composition $b_r,\hdots,b_1$. The composition $\phi(\alpha)$  of $n+3$ is obtained from $L$ after going through the following process:
\begin{itemize}
\item if $b_{r-1}$ is even, then add $1$ to $b_r$; otherwise, append $1$ 
at the beginning of $L$;
\item if $b_1$ is even, then add $2$ to $b_1$; otherwise, append $2$ at the end of $L$.
\end{itemize}

For instance, if $\alpha=UD_2UUDUD_2UDUDUUD$ then we have $B_1=U$, $B_2=D_2U$, $B_3=UDU$, $B_4=D_2UDUDU$, $B_5=UD$, $b_1=1, b_2=2, b_3=3, b_4=6, b_5=2$, $L=2,6,3,2,1$, and $\phi(\alpha) =3,6,3,2,1,2$.

In the case $r=1$, $\alpha$ is either $(UD)^{n/2}$ or $(DU)^{n/2}$. So, we define $\phi((UD)^{n/2})=n+1,2$, and $\phi((DU)^{n/2})=1,n+2$ (these are valid compositions since $n$ has to be even in these cases).

In the case $r=0$,  $\alpha$ is empty, and  we define $\phi(\varepsilon)=1,2$.

\bigskip

Due to the definition, it is clear that the composition $\phi(\alpha)$  belongs to $\cc'(n+3)$, for any $n\geq 0$.

\begin{figure}[h]
\begin{center}
\scalebox{0.8}{\begin{tikzpicture}[ultra thick]
    \fill[bleuclair] (0,-0.4) rectangle (0.4,0.4);
    \fill[rougeclair] (0.4,0.4) rectangle (1.2,-0.4);
    \fill[bleuclair] (1.2,-0.4) rectangle (2.4,0.4);
    \fill[rougeclair] (2.4,0.4) rectangle (4.8,-0.4);
    \fill[bleuclair] (4.8,-0.4) rectangle (5.6,0.4);
    \draw[black, thick] (0,0)--(6.2,0);\draw[black, thick] (0,-0.8)--(0,0.8);
   \draw[black, line width=2pt](0,0)--(0.4,0.4)--(0.8,-0.4)--(1.2,0)--(1.6,0.4)--(2,0)--(2.4,0.4)--(2.8,-0.4)--(3.2,-0)--(3.6,-0.4)--(4,0)--(4.4,-0.4)--(4.8,0)--(5.2,0.4)--(5.6,0);
  % \draw  (1.2,0.8) node {\large $C_2$};\draw  (4.8,0.8) node {\large $C_2$};

\node at (0.2,0.8) {{\color{bleufonce}{$1$}}};
\node at (0.8,0.8) {{\color{rougefonce}{$2$}}};
\node at (1.8,0.8) {{\color{bleufonce}{$3$}}};
\node at (3.6,0.8) {{\color{rougefonce}{$6$}}};
\node at (5.2,0.8) {{\color{bleufonce}{$2$}}};

\draw[-to] (6.7,0)--(7.7,0);

\node at (8.9,0) {$2,6,3,2,1$};

\draw[-to] (10.1,0)--(11.1,0);

\node at (12.5,0) {$3,6,3,2,1,2$};

   \tikzset{every node/.style={circle, draw,fill=white,inner sep=1.5pt}}
   \node at (0,0) {};
   \node at (0.4,0.4) {};
   \node at (0.8,-0.4) {};
   \node at (1.2,0) {};
   \node at (1.6,0.4) {};
   \node at (2,0) {};
   \node at (2.4,0.4) {};
   \node at (2.8,-0.4) {};
   \node at (3.2,-0) {};
   \node at (3.6,-0.4) {};
   \node at (4,0) {};
   \node at (4.4,-0.4) {};
   \node at (4.8,0) {};
    \node at (5.2,0.4) {};
    \node at (5.6,0) {};
\end{tikzpicture}}
\end{center}
\caption{The image by $\phi$ of $\alpha=UD_2UUDUD_2UDUDUUD$ is $\phi(\alpha)=3,6,3,2,1,2$.}
\label{exemple_psi}
\end{figure}

\begin{thm} The map $\phi$ from $\g_n^{[-1,1]}$  to  $\cc'(n+3)$ is a bijection.
\end{thm}

\begin{proof}
Since $\g_n^{[-1,1]}$ and $\cc'(n+3)$ have the same cardinality (see the end of subsection \ref{enumeration_-t_t}), it suffices to prove that $\phi$ is surjective.

 Let $c=c_1,\hdots,c_r$ be an element of $\cc'(n+3)$, with $r\geq 2$ (the case $r=1$ does not occur since $c_1$ is odd and $c_r$ is even implies $r\geq 2$). Thus, $r$ is necessarily even. We distinguish four cases: ($i$) $c_1=1$ and $c_r=2$, ($ii$) $c_1=1$ and $c_r$ is even and greater than $2$,  ($iii$) $c_1$ is odd and greater than $1$ and $c_r=2$, ($iv$) $c_1$ is odd and greater than $1$ and $c_r$ is even and greater than $2$.

 According to each case, we define $\alpha\in\g_n^{[-1,1]}$ such that $\phi(\alpha)=c$:

\noindent Case ($i$): Since  $c_{r-1}$ is odd, $c_{r-2}$ is even, and so on, and finally $c_2$ is even. 
$$\alpha=
\colorbox{bleuclair}{$(UD)^{(c_{r-1}-1)/2}U$} \colorbox{rougeclair}{$D_2U(DU)^{(c_{r-2}-2)/2}$} \hdots \colorbox{bleuclair}{$(UD)^{(c_3-1)/2}U$} \colorbox{rougeclair}{$D_2U(DU)^{(c_2-2)/2}$}.
$$
%This GDAP has indeed length $(n+3)-c_1-c_r=n$, and its image by $\phi$ is %$c=1,c_2,\hdots,c_{r-1},2$. 

\noindent Case ($ii$): 
$$\alpha=
\colorbox{rougeclair}{$(DU)^{(c_r-2)/2}$}
\colorbox{bleuclair}{$(UD)^{(c_{r-1}-1)/2}U$} \hdots \colorbox{bleuclair}{$(UD)^{(c_3-1)/2}U$} \colorbox{rougeclair}{$D_2U(DU)^{(c_2-2)/2}$}.
$$
%Once again, this GDAP has indeed length $(n+3)-c_1-2=n$, and its image by $\phi$ is %$c=1,c_2,\hdots,(c_r-2)+2$.
\noindent Case ($iii$):
$$\alpha=
\colorbox{bleuclair}{$(UD)^{(c_{r-1}-1)/2}U$}
\colorbox{rougeclair}{$D_2U(DU)^{(c_{r-2}-2)/2}$} \hdots \colorbox{rougeclair}{$D_2U(DU)^{(c_2-2)/2}$}
\colorbox{bleuclair}{$(UD)^{(c_1-1)/2}$}.
$$

\noindent Case ($iv$):
$$\alpha=
\colorbox{rougeclair}{$(DU)^{(c_r-2)/2}$}
\colorbox{bleuclair}{$(UD)^{(c_{r-1}-1)/2}U$} \hdots \colorbox{rougeclair}{$D_2U(DU)^{(c_2-2)/2}$}
\colorbox{bleuclair}{$(UD)^{(c_1-1)/2}$}.
$$

For each case, it is clear that $\alpha$ belongs to $\g_n^{[-1,1]}$, which implies that $\phi$ is surjective, and then bijective.
\end{proof}

\section{DAP with a special first return decomposition}

Recently in \cite{bakipe},  the authors introduced and enumerated the subset $\mathcal{D}^{h,\geq}$ of restricted Dyck paths defined as follows: the set $\mathcal{D}^{h,\geq}$  is the union of the empty Dyck path with all Dyck paths $P$ having a first return decomposition $P=U\alpha D\beta$ satisfying the  conditions:
\begin{equation}\left\{\begin{array}{l}
\alpha,\beta\in \mathcal{D}^{h,\geq},\\
h(U\alpha D)\geq h(\beta),\\
\end{array}\right.
\label{eq1}
\end{equation}
 where $h(\alpha)$ is the maximal ordinate reached by the path $\alpha$.
The authors prove algebraically and bijectively  that $n$-length paths in $\mathcal{D}^{h,\geq}$ are in one-to-one correspondence with  Motzkin paths of length $n$. Based on this decomposition and in the same way as for Dyck paths, we define a subset of $\mathcal{A}\cup\{\varepsilon\}\subset\mathcal{G}$ as follows.  
 The set $\mathcal{H}$  is the union of the empty path with all DAP $\gamma\in\mathcal{A}$ having a first return decomposition satisfying the following condition:

   ($C$) $\gamma= \alpha\beta \mbox{ with } \alpha\in\mathcal{P}\cup\{UD\}, \mbox{ and }\alpha^\flat\in\mathcal{H} \mbox{ whenever } \alpha\neq UD,  \beta\in \mathcal{H} \mbox{ and } h(\alpha )\geq h(\beta).$
   
  For $n\geq 0$, we denote by $\mathcal{H}_n$ the set of DAP  of length $n$ in $\mathcal{H}$.  For instance, we have $\mathcal{H}_0=\{\varepsilon\}$, $\mathcal{H}_1=\emptyset$, $\mathcal{H}_2=\{UD\}$, $\mathcal{H}_3=\{UUD_2\}$, $\mathcal{H}_4=\{UUUD_3, UDUD\}$,
  $\mathcal{H}_5=\{UUUUD_4, UUDUD_2, UUD_2UD\}$.

In this section, we enumerate the set $\mathcal{H}_n$.  For $k\geq 0$, let $A_k(x)=\sum_{n\geq 0} a_{n,k}x^n$ (resp. $B_k(x)=\sum_{n\geq 0} b_{n,k}x^n$) be the generating function where the coefficient $a_{n,k}$ (resp. $b_{n,k}$) is the number of DAP in $\mathcal{H}_n$ having a maximal height equal to $k$ (resp. of at most $k$). So, we have $B_k(x)=\sum\limits_{i=0}^kA_i(x)$ and the generating function for the set $\mathcal{H}$, namely $B(x)$, is given by $B(x)=\lim\limits_{k\rightarrow\infty}B_k(x)$.

Due to the definition of $\mathcal{H}$,  we have
$$\left\{\begin{array}{ll}
A_0(x)&=B_0(x)=1,\\
A_1(x)&=x^2A_0(x)B_1(x),\\
A_k(x)&=xA_{k-1}(x)B_k(x).
\end{array}\right.$$

\begin{lemme} For $k\geq 1$, we have
$$B_{k-1}(x)=\frac{(1-x^3+x)B_k(x)-1}{x^2B_k(x)+x}.$$
\end{lemme}
\begin{proof}
We proceed by induction on $k$. Since $B_0(x)=1$ and $B_1(x)=\frac{x^2}{1-x^2}$ it is easy to check that $B_{0}(x)=\frac{(1-x^3+x)B_1(x)-1}{x^2B_1(x)+x}$.

Now, assume that $B_{i-1}(x)=\frac{(1-x^3+x)B_i(x)-1}{x^2B_i(x)+x}$ for $1\leq i \leq k-2$, we prove the result for $i=k-1$. From the above equations, and the recurrence hypothesis on  $B_{k-2}(x)$, we obtain $$\begin{array}{ll}
B_k(x)&=A_k(x)+B_{k-1}(x)\\
      &=xA_{k-1}(x)B_k(x)+B_{k-1}(x)\\
      &=x(B_{k-1}-B_{k-2})B_k(x)+B_{k-1}(x)\\
      &=x(B_{k-1}-\frac{(1-x^3+x)B_{k-1}(x)-1}{x^2B_{k-1}(x)+x})B_k(x)+B_{k-1}(x).\\
\end{array}$$
Isolating $B_{k-1}(x)$, we obtain $B_{k-1}(x)=\frac{(1-x^3+x)B_k(x)-1}{x^2B_k(x)+x}$, which completes the induction.
\end{proof}

Taking the limit in the relation of Lemma 1 whenever $k$ tends to $\infty$,  we obtain 
$$B(x)=\frac{(1+x-x^3)B(x)-1}{x^2B(x)+x},$$ which induces the following result.
\begin{thm} The o.g.f. that counts the set $\mathcal{H}$ with respect to the length is given by $$B(x)={\frac {1-{x}^{3}-\sqrt {{x}^{6}-2\,{x}^{3}-4\,{x}^{2}+1}}{2{x}^{
2}}}.
$$
\end{thm}

The above generating function $B(x)$ counts also Motzkin paths of length $n$ avoiding the patterns $UH$, $HU$ and $HH$, see sequence 
\oeis{A329699}. The first terms of its series expansion are:
$1+x^2+x^3+2x^4+3x^5+6x^6+10x^7+20x^8+36x^9+72x^{10}+136x^{11}+273x^{12}$.

We finish this part with a natural question.

\noindent 
{\bf Open question:} 
The sets $\mathcal{H}$ and that of 
$\{UH, HU, HH\}$-avoiding Motzkin paths are thus in bijection and
it would be interesting to exhibit a constructive bijection between them.

\end{document}